\documentclass[11pt]{article}

\usepackage{amssymb,latexsym}
\usepackage{epsfig}
\usepackage{eufrak}
\usepackage{amsmath}
\usepackage{mathrsfs}
\usepackage{color}

\newtheorem{theorem}{Theorem}
\newtheorem{lemma}{Lemma}

\newcommand{\be}{\begin{equation}}
\newcommand{\ee}{\end{equation}}
\newcommand{\bee}{\begin{eqnarray*}}
\newcommand{\eee}{\end{eqnarray*}}
\newcommand{\bel}{\begin{eqnarray}}
\newcommand{\eel}{\end{eqnarray}}
\newcommand{\bec}{\begin{cases}}
\newcommand{\eec}{\end{cases}}
\newcommand{\bem}{\begin{bmatrix}}
\newcommand{\eem}{\end{bmatrix}}

\newcommand{\la}{\label}
\newcommand{\li}{\left}
\newcommand{\ri}{\right}

\newcommand{\ovl}{\overline}
\newcommand{\udl}{\underline}

\newcommand{\lc}{\lceil}
\newcommand{\rc}{\rceil}
\newcommand{\lf}{\lfloor}
\newcommand{\rf}{\rfloor}

\newcommand{\vep}{\varepsilon}

\newcommand{\de}{\delta}

\newcommand{\ga}{\gamma}

\newcommand{\se}{\theta}

\newcommand{\ze}{\zeta}

\newcommand{\vro}{\varrho}

\newcommand{\om}{\omega}
\newcommand{\Om}{\Omega}

\newcommand{\f}{\frac}

\newcommand{\cd}{\cdots}

\newcommand{\qu}{\quad}
\newcommand{\qqu}{\qquad}

\newcommand{\mscr}{\mathscr}

\newcommand{\bb}{\mathbb}

\newcommand{\bs}{\boldsymbol}

\newcommand{\tx}{\text}

\newcommand{\iy}{\infty}

\newcommand{\pa}{\partial}

\newcommand{\bed}{\begin{description}}
\newcommand{\eed}{\end{description}}
\newcommand{\bei}{\begin{itemize}}
\newcommand{\eei}{\end{itemize}}
\newcommand{\ben}{\begin{enumerate}}
\newcommand{\een}{\end{enumerate}}
\newcommand{\bib}{\bibitem}
\newcommand{\beL}{\begin{lemma}}
\newcommand{\eeL}{\end{lemma}}
\newcommand{\beT}{\begin{theorem}}
\newcommand{\eeT}{\end{theorem}}
\newcommand{\sect}{\section}

\newcommand{\bpf}{\begin{pf}}
\newcommand{\epf}{\end{pf}}
\newcommand{\bsk}{\bigskip}

\newcommand{\pfbox}{\hfill\mbox{$\Box$}}

\newenvironment{pf}{\paragraph*{Proof{\rm.}}}{\pfbox\bigskip}

\begin{document}

\title{{\bf Interval Estimation of Bounded Variable Means via Inverse Sampling}
\thanks{The author is currently with Department of Electrical Engineering,
Louisiana State University at Baton Rouge, LA 70803, USA, and Department of Electrical Engineering, Southern
University and A\&M College, Baton Rouge, LA 70813, USA; Email: chenxinjia@gmail.com}}

\author{Xinjia Chen}

\date{February, 2008}

\maketitle

\begin{abstract}

In this paper, we develop interval estimation methods for means of bounded random variables based on a sequential procedure such that the
sampling is continued until  the sample sum is no less than a prescribed threshold.

\end{abstract}

\sect{Inverse Sampling}

It is a ubiquitous problem to estimate the means of bounded random variables.  Specially,  the problem of estimating the probability of an event
can be formulated as the estimation of the mean of a Bernoulli  variable.  Moreover, in many applications, one needs to estimate a quantity
$\mu$ which can be bounded in $[0,1]$ after proper operations of scaling and translation.  A typical approach is to design an experiment that
produces a random variable $X$ distributed in $[0, 1]$ with expectation $\mu$, run the experiment independently a number of times, and use the
average of the outcomes as the estimate \cite{Dagum}.

The objective of this paper is to develop interval estimation methods for means of bounded variables based a sequential sampling scheme
described as follows.

Let $X \in [0, 1]$ be a bounded random variable defined in a probability space $(\Om, \mscr{F}, \Pr)$ with mean value  $\bb{E}[X] = \mu \in (0,
1)$. We wish to estimate the mean of $X$ by using a sequence of i.i.d. random samples $X_1, \; X_2, \; \cd$ of $X$ based on the following {\it
inverse sampling} scheme:

{\it Continue sampling until the sample size reach a number $\bs{n}$ such that the sample sum $\sum_{i = 1}^{\bs{n}} X_i$ is no less than a
positive number $\ga$}.

\bsk

We call this an {\it inverse sampling} scheme, since it reduces to the classical {\it inverse binomial sampling} scheme \cite{H} \cite{H2} in
the special case that $X$ is a Bernoulli random variable.

\sect{Interval Estimation}

When the experiment of inverse sampling is completed, it is desirable to construct a confidence interval for the mean value $\mu$ in terms of
the random sampling number $\bs{n}$.  For this purpose, we have

\beT \la{thm1}
Let $\mscr{H} (z, \mu) =  \ln \li ( \f{\mu}{z} \ri ) + \li ( \f{1}{z} - 1 \ri ) \ln \li (  \f{ 1 - \mu } { 1 - z } \ri )$ for $0 <
z < 1$ and $0 < \mu < 1$.  Let $\de \in (0,1)$.  Let $\ovl{\bs{\mu}}$ be a random variable such that $\ovl{\bs{\mu}} = 1$ for $\bs{n} \leq \ga +
1$ and that
\[
\mscr{H} \li ( \f{ \ga}{ \bs{n} - 1},  \ovl{\bs{\mu}} \ri ) = \f{ \ln \f{\de}{2} } {  \ga }, \qqu \f{ \ga}{ \bs{n} - 1} < \ovl{\bs{\mu}} < 1
\]
for $\bs{n} > \ga + 1$.   Let $\udl{\bs{\mu}}$ be a random variable such that
\[
\mscr{H} \li ( \f{ \ga}{ \bs{n} },  \udl{\bs{\mu}} \ri ) = \f{ \ln \f{\de}{2} } {  \ga }, \qqu 0 < \udl{\bs{\mu}} < \f{ \ga}{ \bs{n}}.
\]
Then, \[ \Pr \{  \udl{\bs{\mu}} < \mu < \ovl{\bs{\mu}}  \} \geq 1 - \de.
\]

 \eeT

The proof is given in Appendix A.  It should be noted that, due to the monotone property of the function $\mscr{H} (z, \mu)$ with respect to
$\mu$, the confidence limits can be readily determined by a bisection search method.

In the special case that $X$ is a Bernoulli random variable and $\ga$ is an integer, the confidence interval can be constructed in a slightly
different way as follows.

 \beT
\la{thm2} Let $\de \in (0,1)$.  Let $\ovl{\bs{\mu}}$ be a random variable such that $\ovl{\bs{\mu}} = 1$ for $\bs{n} = \ga$ and that
\[
\mscr{H} \li ( \f{ \ga}{ \bs{n} },  \ovl{\bs{\mu}} \ri ) = \f{ \ln \f{\de}{2} } {  \ga }, \qqu \f{ \ga}{ \bs{n} } < \ovl{\bs{\mu}} < 1
\]
for $\bs{n} > \ga$.  Let $\udl{\bs{\mu}}$ be a random variable such that
\[
\mscr{H} \li ( \f{ \ga}{ \bs{n} },  \udl{\bs{\mu}} \ri ) = \f{ \ln \f{\de}{2} } {  \ga }, \qqu 0 < \udl{\bs{\mu}} < \f{ \ga}{ \bs{n}}.
\]
Then, \[ \Pr \{  \udl{\bs{\mu}} < \mu < \ovl{\bs{\mu}}  \} \geq 1 - \de.
\]

 \eeT

The proof is given in Appendix B.  As mentioned earlier, the confidence limits can be readily determined by a bisection search method.

Theorems \ref{thm1}  and \ref{thm2} are established by employing Hoeffding's inequality \cite{Hoeffding}. If we replace the Hoeffding's
inequality by Massart's inequality (i.e., Theorem 2 at page 1271 of \cite{Massart:90}), which is slightly more conservative, we can obtain via
analogy arguments explicit formulas for interval estimation.  In this regard, we have Theorem \ref{thm3} for the general inverse sampling
scheme.

\begin{theorem}
\la{thm3} Let $\delta \in (0,1)$ and $\theta = \frac{9}{ 2 \ln \frac{2}{\delta} }$. Define \bee \udl{\bs{\mu}} & = & \f{\ga}{\bs{n}} +
\frac{3}{4 + \bs{n} \theta}
\left[ 1 -  \f{2 \ga}{\bs{n}} - \sqrt{ 1 + \se \ga \li (1 - \f{\ga}{\bs{n}} \ri ) } \right],\\
\ovl{\bs{\mu}} & = & \bec \f{\ga}{\bs{n} - 1} + \frac{3}{4 + (\bs{n} - 1) \theta} \left[ 1 -  \f{2 \ga}{\bs{n} - 1} + \sqrt{ 1 + \se \ga \li (1
- \f{\ga}{\bs{n} - 1} \ri ) } \right] & \tx{for} \; \bs{n} > \ga + 1,\\
1 & \tx{for} \; \bs{n} \leq \ga + 1. \eec \eee Then,
\[
\Pr \{ \udl{\bs{\mu}} < \mu  < \ovl{\bs{\mu}} \} \geq 1 - \delta.
\]
\end{theorem}

For the inverse binomial sampling scheme (with integer $\ga$), we have

\begin{theorem}
\la{thm4} Let $\delta \in (0,1)$ and $\theta = \frac{9}{ 2 \ln \frac{2}{\delta} }$. Define
\[
\udl{\bs{\mu}} = \f{\ga}{\bs{n}} + \frac{3}{4 + \bs{n} \theta} \left[ 1 -  \f{2 \ga}{\bs{n}} - \sqrt{ 1 + \se \ga \li (1 - \f{\ga}{\bs{n}} \ri )
} \right],
\]
\[
\ovl{\bs{\mu}} = \f{\ga}{\bs{n}} + \frac{3}{4 + \bs{n} \theta} \left[ 1 -  \f{2 \ga}{\bs{n}} + \sqrt{ 1 + \se \ga \li (1 - \f{\ga}{\bs{n}} \ri )
} \right].
\]
Then,
\[
\Pr \{ \udl{\bs{\mu}} < \mu < \ovl{\bs{\mu}} \} \geq 1 - \delta.
\]
\end{theorem}

\sect{Conclusion}

We have established rigorous  and simple interval estimation methods for means of bounded random variables.  The construction of confidence
intervals is based on inverse sampling and requires little computation. The nominal coverage probability of confidence intervals  is always
guaranteed.

\appendix

 \sect{Proof Theorem 1}

We need some preliminary results.

The following lemma is a classical result, known as Hoeffding's inequality \cite{Hoeffding}.

 \beL \la{lem1}  Let $X_1, \cd, X_n$ be i.i.d. random variables
bounded in $[0,1]$ with common mean value $\mu \in (0,1)$. Then, {\small $\Pr \li \{ \f{\sum_{i =1}^n X_i}{n} \geq z \ri \} \leq \exp \li ( n z
\mscr{H} (z, \mu) \ri )$ } for $1
> z > \mu = \bb{E}[ X_i]$. Similarly, {\small $\Pr \li \{  \f{\sum_{i =1}^n X_i}{n} \leq z  \ri \} \leq \exp \li
(  n z \mscr{H} (z, \mu) \ri )$ } for $0 < z < \mu = \bb{E}[ X_i]$. \eeL

The following lemma has been established by Chen \cite{Chen}.

 \beL \la{lem2} $\mscr{H}( (1 + \vep) \mu, \mu)$ is monotonically decreasing with respect
to $\vep \in \li ( 0, \; \f{1}{\mu} - 1 \ri )$. Similarly, $\mscr{H}( (1 - \vep) \mu, \mu)$ is monotonically decreasing with respect to $\vep
\in \li ( 0, \; 1 \ri )$.
 \eeL

\beL \la{lem3} For any $\ga > 0 $ and $\vep \in (0,1)$,
\[
\Pr \li \{  \bs{n} \leq \f{ \ga } { \mu (1 + \vep) } \ri \}  \leq \exp( \ga \mscr{H} ((1 + \vep) \mu, \mu) ).
\]
\eeL

\bpf  Since $\bs{n}$ is an integer, we have
\[
 \Pr \li \{  \bs{n} \leq \f{ \ga } { \mu (1 +
\vep) } \ri \} = \Pr \li \{  \bs{n} \leq \li \lf \f{ \ga } { \mu (1 + \vep) } \ri \rf \ri \} = \Pr \li \{ \bs{n}
\leq \f{ \ga } { \mu (1 + \vep^*) } \ri \}
\]
where $\vep^*$ is a number depends on $\mu$ and $\vep$ such that {\small $\f{ \ga } { \mu (1 + \vep^*) } = \li \lf \f{ \ga } { \mu (1 + \vep) }
\ri \rf$}.  Here $\lf . \rf$ denotes the floor function.   Clearly, {\small $\vep^* = \f{ \ga } { \mu  \li \lf \f{ \ga } { \mu (1 + \vep) } \ri
\rf} - 1 \geq \vep
> 0$}. For simplicity of notation, let {\small $m = \f{ \ga } { \mu (1 + \vep^*) }$}. Since $m$ is a nonnegative integer, it can be zero or a
natural number.  If $m = 0$, then
\[
\Pr \li \{  \bs{n} \leq \f{ \ga } { \mu (1 + \vep) } \ri \}  = \Pr \{ \bs{n} \leq m \} = 0 < \exp( \ga \mscr{H} ((1 + \vep) \mu, \mu) ).
\]
Otherwise if $m \geq 1$, then
\[
\Pr \li \{ \bs{n} \leq \f{ \ga } { \mu (1 + \vep^*) } \ri \} = \Pr \{ \bs{n} \leq m \} = \Pr \{ X_1 + \cd + X_m
\geq \ga \} = \Pr \{ \ovl{X} \geq z \},
\]
where {\small $\ovl{X} = \f{ \sum_{i=1}^m X_i } { m }$} and {\small $z = \f{\ga}{ m } = \mu (1 + \vep^*) >
\mu$}. Now we shall consider three cases.

\bed

\item  (i): In the case of $z > 1$, we have {\small $\Pr \{ \ovl{X} \geq z \} \leq \Pr \li \{ \sum_{i = 1}^m X_i > m
\ri \} = 0 < \exp( \ga \mscr{H} ((1 + \vep) \mu, \mu) )$}.

\item (ii): In the case of $z = 1$, we have $\mu = \f{1}{1 + \vep^*}, \; m = \ga$ and {\small \bee \Pr \{ \ovl{X} \geq z
\} & = & \Pr \li \{ \sum_{i = 1}^m X_i = m \ri \} = \prod_{i = 1}^m \Pr \{ X_i
= 1 \} \leq  \prod_{i = 1}^m \bb{E}[ X_i ] = \mu^m\\
& = & \li ( \f{ 1 } {1 + \vep^* }  \ri )^\ga \leq \li ( \f{ 1 } {1 + \vep}  \ri )^\ga \leq \exp( \ga \mscr{H} ((1 + \vep) \mu, \mu) ). \eee}

\item (iii): In the case of  $\mu < z < 1$, by Lemma \ref{lem1}, we have
\[
\Pr \{ \ovl{X} \geq z \} \leq \exp( m z \mscr{H} (z, \mu) ) = \exp( \ga \mscr{H} ((1 + \vep^*) \mu, \mu) ).
\]
Since $\vep^* \geq \vep$,  it must be true that $\mu (1 + \vep) \leq \mu (1 + \vep^*) < 1$ and that $\mscr{H} ((1 + \vep^*) \mu, \mu) \leq
\mscr{H} ((1 + \vep) \mu, \mu)$ as a result of Lemma \ref{lem2}. Hence, \be \la{remm} \Pr \li \{  \bs{n} \leq \f{ \ga } { \mu (1 + \vep) } \ri
\} = \Pr \{ \ovl{X} \geq z \} \leq \exp( \ga \mscr{H} ((1 + \vep) \mu, \mu) ). \ee

\eed

 Therefore, we have shown
 \be
 \la{taila}
\Pr \li \{  \bs{n} \leq \f{ \ga } { \mu (1 + \vep) } \ri \}  \leq \exp( \ga \mscr{H} ((1 + \vep) \mu, \mu) )
 \ee
 for all cases.

 \epf

\bsk

\beL  \la{lem4} For any $\ga > 0 $ and $\vep \in (0,1)$,
\[
\Pr \li \{  \bs{n} \geq \f{ \ga } { \mu (1 - \vep) } \ri \} \leq \exp \li ( \ga \mscr{H} \li ( \f{\ga}{ \f{ \ga } { \mu (1 -\vep) } - 1 }, \mu
\ri )  \ri ).
\]
\eeL

\bpf

Since $\bs{n}$ is an integer, we have
\[
\Pr \li \{  \bs{n} \geq \f{ \ga } { \mu (1 - \vep) } \ri \} =
 \Pr \li \{  \bs{n} \geq \li \lc \f{ \ga } { \mu (1
-\vep) } \ri \rc \ri \} = \Pr \li \{  \bs{n} > \li \lc \f{ \ga } { \mu (1 -\vep) } - 1 \ri \rc \ri \}
\]
where $\lc . \rc$ denotes the ceiling function.  Let $\ze$ be a number such that  {\small $\f{ \ga } { \mu (1 - \vep) } - 1 = \f{ \ga } { \mu (1
- \ze) }$} for any $\mu \in (0,1)$. Hence,
\[
\Pr \li \{  \bs{n} > \li \lc \f{ \ga } { \mu (1 -\vep) } - 1 \ri \rc \ri \} = \Pr \li \{  \bs{n} > \li \lc \f{
\ga } { \mu (1 - \ze ) }  \ri \rc \ri \} = \Pr \li \{  \bs{n} >  \f{ \ga } { \mu (1 - \ze^* ) }  \ri \}
\]
with $\ze^*$ satisfying {\small $\f{ \ga } { \mu (1 - \ze^* ) }  = \li \lc \f{ \ga } { \mu (1 - \ze ) } \ri
\rc$}. Clearly, $1 > \ze^* \geq \ze > 0$.  Let {\small $m = \f{ \ga } { \mu (1 - \ze^* ) }$}. Then, $m$ is a
positive integer and
\[
\Pr \{  \bs{n} > m \} = \Pr \{ X_1 + \cd + X_m < \ga \} = \Pr \{ \ovl{X} < z \}
\]
where {\small $\ovl{X} = \f{ \sum_{i=1}^m X_i } { m }$} and $z = (1 - \ze^*) \mu$.  Applying Lemma \ref{lem1}, we have {\small \bee \Pr \li \{
\bs{n} > \li \lc \f{ \ga }
 { \mu (1 - \ze) } \ri \rc  \ri \}  & = & \Pr \{ \ovl{X} < z \} \leq  \exp ( m z
\mscr{H} ( z, \mu) )\\
&  = & \exp ( \ga \mscr{H} ( (1 - \ze^*) \mu, \mu) ). \eee} Note that $\mscr{H} ( (1 - \ze^*) \mu, \mu) \leq \mscr{H} ( (1 - \ze) \mu, \mu)$ as
a result of $1 > \ze^* \geq \ze > 0$ and Lemma \ref{lem2}. Hence,
\[
\Pr \{ \ovl{X} < z \} \leq \exp ( \ga \mscr{H} ( (1 - \ze) \mu, \mu) ) = \exp \li ( \ga \mscr{H} \li ( \f{\ga}{ \f{ \ga } { \mu (1 -\vep) } - 1
}, \mu \ri )  \ri ).
\]

\epf

\beL \la{lem5} For any $\ga > 0 $,  {\small $\mscr{H} \li ( \f{\ga}{ \f{ \ga } { \mu (1 -\vep) } - 1 }, \mu \ri )$} is monotonically decreasing
with respect to {\small $\vep \in \li ( \f{ \mu } { \ga + \mu}, 1 \ri ) $}. \eeL

\bpf

Let \[ \vro = \f{1}{1 - \vep} - 1.
\]
Then, \[ \mscr{H} \li ( \f{\ga}{ \f{ \ga } { \mu (1 -\vep) } - 1 }, \mu \ri ) = \mscr{H} \li ( \f{\ga \mu}{\ga (1 + \vro) - \mu}, \mu \ri ).
\]
 Let {\small $z = \f{\ga \mu}{\ga (1 + \vro)  - \mu}$} and $m = \f{\ga}{z}$.  For $\vro > \f{ \mu } { \ga }$, we have $0 < z < \mu, \; \f{ \pa
m } { \pa \vro } > 0$ and  {\small \bee \f{ \pa [\ga \mscr{H}(z, \mu)] } { \pa \vro } & = & \ga \li [- \f{1}{z} \f{ \pa z } { \pa \vro } + \li
(\f{1}{z} - 1 \ri ) \f{1}{1 - z} \f{ \pa z } { \pa \vro } - \f{1}{z^2} \f{ \pa z } {
\pa \vro } \ln \li ( \f{1 - \mu}{1 - z} \ri ) \ri ]\\
& = & - \f{\ga}{z^2} \f{ \pa z } { \pa \vro } \ln \li ( \f{1 - \mu}{1 - z} \ri ) =  - \f{m \ga}{m z^2} \f{ \pa z } { \pa \vro } \ln \li ( \f{1 - \mu}{1 - z} \ri )\\
& = & -\f{m }{ z } \f{ \pa z } { \pa \vro } \ln \li ( \f{1 - \mu}{1 - z} \ri ) =  \f{ \pa m } { \pa \vro } \ln \li ( \f{1 - \mu}{1 - z} \ri ) <
0, \eee} which implies that $\mscr{H}(z, \mu)$ is monotonically decreasing with respect to {\small $\vro > \f{ \mu } { \ga }$}.  By the relation
between $\vro$ and $\vep$, we have that $\vro$ increases as $\vep$ increases and that $\vro > \f{ \mu } { \ga }$ if and only if $\vep > \f{ \mu
} { \ga + \mu}$.  This proves the lemma.

\epf

The following lemma can be shown by direct computation.

 \beL \la{lem6}
\[
\f{\pa \mscr{H} (z, \mu)  } {\pa \mu} = \f{1}{\mu (1 - \mu) } - \f{1}{z (1 - \mu)},
\] which is negative for $1 > \mu > z > 0$, and positive for $0 < \mu < z < 1$.  Moreover,
\[
\lim_{\mu \to z} \mscr{H} (z, \mu) = 0, \qqu \lim_{\mu \to 0} \mscr{H} (z, \mu)  = - \iy, \qqu \lim_{\mu \to 1} \mscr{H} (z, \mu)  = - \iy
\]
for $0 < z < 1$.  \eeL

\beL \la{lem7} For any $\ga > 0 $ and $\de \in (0,1)$,  there exists $\vep^* \in \li ( \f{\mu}{\ga + \mu},  1 \ri )$ such that
\[
\mscr{H} \li ( \f{ \ga}{ \f{ \ga } { \mu (1 - \vep^*)  } - 1},  \mu \ri ) = \f{ \ln \f{\de}{2} } {  \ga }.
\]
\eeL

\bpf

Note that
\[
\lim_{\vep \to \f{\mu}{\ga + \mu} } \mscr{H} \li ( \f{ \ga}{ \f{ \ga } { \mu (1 - \vep)  } - 1},  \mu \ri ) = \mscr{H} \li ( \mu,  \mu \ri ) =
0, \qu \lim_{\vep \to 1 } \mscr{H} \li ( \f{ \ga}{ \f{ \ga } { \mu (1 - \vep)  } - 1},  \mu \ri ) = - \iy.
\]
By Lemma \ref{lem5}, {\small $\mscr{H} \li ( \f{ \ga}{ \f{ \ga } { \mu (1 - \vep)  } - 1},  \mu \ri )$} is monotonically decreasing with respect
to $\vep \in \li (\f{\mu}{\ga + \mu},  1 \ri )$. Since $- \iy < \f{ \ln \f{\de}{2} } {  \ga } < 0$, the existence of $\vep^*$ is established.

\epf

\beL \la{lem8} For any $\ga > 0 $ and $\de \in (2 \mu^\ga,1)$,  there exists $\vep^* \in \li (0, \f{1}{\mu} - 1 \ri )$ such that
\[
\mscr{H} ( (1 + \vep^*) \mu, \mu ) = \f{ \ln \f{\de}{2} } {  \ga }.
\]
\eeL

\bpf

Note that
\[
\lim_{\vep \to 0} \mscr{H} ( (1 + \vep) \mu, \mu ) = 0, \qu  \lim_{\vep \to \f{1}{\mu} - 1  } \mscr{H} ( (1 + \vep) \mu, \mu ) = \ln \mu.
\]
By Lemma \ref{lem2}, {\small $\mscr{H} ( (1 + \vep) \mu, \mu )$} is monotonically decreasing with respect to $\vep \in \li (0, \f{1}{\mu} - 1
\ri ) $. Since $\ln \mu < \f{ \ln \f{\de}{2} } {  \ga } < 0$ for $\de \in (2 \mu^\ga,1)$, the existence of $\vep^*$ is established.

\epf

\beL \la{lem9} For any $\ga > 0 $ and $\de \in (0,1)$,
\[
\Pr \li \{  \mscr{H} \li ( \f{ \ga}{ \bs{n} - 1},  \mu \ri ) \leq \f{ \ln \f{\de}{2} } {  \ga }, \; \; \mu \geq \f{ \ga}{ \bs{n} - 1}, \; \bs{n}
> \ga + 1 \ri \} \leq \f{\de}{2}.
\]
\eeL

\bpf

By Lemma \ref{lem7},  there exists $\vep^* \in \li ( \f{\mu}{\ga + \mu},  1 \ri )$ such that {\small $\mscr{H} \li ( \f{ \ga}{ \f{ \ga } { \mu
(1 - \vep^*) } - 1},  \mu \ri ) = \f{ \ln \f{\de}{2} } {  \ga }$}.  By Lemma \ref{lem4},
\[
\Pr \li \{ \bs{n} \geq \f{ \ga } { \mu (1 - \vep^*)  } \ri \} \leq \exp \li ( \ga  \mscr{H} \li ( \f{ \ga}{ \f{ \ga } { \mu (1 - \vep^*) } - 1},
\mu \ri )  \ri ) = \f{\de}{2}.
\]
Therefore, to show Lemma \ref{lem9}, it suffices to show
\[
\li \{  \mscr{H} \li ( \f{ \ga}{ \bs{n} - 1},  \mu \ri ) \leq \f{ \ln \f{\de}{2} } {  \ga }, \; \; \mu \geq \f{ \ga}{ \bs{n} - 1}, \; \bs{n}
> \ga + 1 \ri \} \subseteq \li \{ \bs{n} \geq \f{ \ga } { \mu (1 - \vep^*)  } \ri \}.
\]
Let $\om \in \li \{  \mscr{H} \li ( \f{ \ga}{ \bs{n} - 1},  \mu \ri ) \leq \f{ \ln \f{\de}{2} } {  \ga }, \; \; \mu \geq \f{ \ga}{ \bs{n} - 1},
\; \bs{n} > \ga + 1 \ri \}$ and $n = \bs{n}(\om)$.  Let
\[
\vep =  1 - \f{ \ga } { n \mu }.
\]
Note that
\[
\mu \geq \f{ \ga}{ n - 1} > 0,
\]
which can be written as \[ 1 - \f{ \ga } { n \mu } \geq \f{\mu}{\ga + \mu}.
\]
Hence,
\[
1 > \vep \geq \f{\mu}{\ga + \mu}.
\]
By the definition of $\vep$, we have
\[
\f{ \ga}{ n - 1} = \f{ \ga}{ \f{ \ga } { \mu (1 - \vep)  } - 1}
\]
and thus
\[
\mscr{H} \li ( \f{ \ga}{ \f{ \ga } { \mu (1 - \vep)  } - 1},  \mu \ri ) \leq \f{ \ln \f{\de}{2} } {  \ga }.
\]
As a result, $\mscr{H} \li ( \f{ \ga}{ \f{ \ga } { \mu (1 - \vep)  } - 1},  \mu \ri ) \leq \mscr{H} \li ( \f{ \ga}{ \f{ \ga } { \mu (1 - \vep^*)
} - 1},  \mu \ri )$, and by Lemma \ref{lem5}
\[
\vep \geq \vep^*.
\]
By this inequality and  the definition of $\vep$,
\[
n = \f{ \ga } { \mu (1 - \vep)  } \geq \f{ \ga } { \mu (1 - \vep^*)  }.
\]
So,
\[
\om \in \li \{ \bs{n} \geq \f{ \ga } { \mu (1 - \vep^*)  } \ri \}.
\]
This shows the inclusion relationship of the sets. The lemma is thus proved.  \epf

\beL \la{lem10}  For any $\ga > 0 $ and $\de \in (0,1)$,
\[
\Pr \li \{  \mscr{H} \li ( \f{ \ga}{ \bs{n} },  \mu \ri ) \leq \f{ \ln \f{\de}{2} } {  \ga }, \; \; \mu < \f{ \ga}{ \bs{n} } \ri \} \leq
\f{\de}{2}.
\]
\eeL

\bpf

There are two cases: Case (i) $0 < \de < 2 \mu^\ga$; Case (ii) $\de \geq 2 \mu^\ga$.

We first consider Case (i). Note that $\mscr{H}(z, \mu)$ is monotonically decreasing with respect to $z \in (\mu, 1)$ and that
\[
\lim_{z \to 1} \mscr{H}(z, \mu) = \ln \mu.
\]
Since $0 < \de < 2 \mu^\ga$, we have $\ln \mu > \f{ \ln \f{\de}{2} } { \ga }$.  As a result, we have that
\[
\mscr{H}(z, \mu) > \f{ \ln \f{\de}{2} } { \ga }
\]
for any $z \in (\mu, 1)$.  This implies that $\li \{  \mscr{H} \li ( \f{ \ga}{ \bs{n} },  \mu \ri ) \leq \f{ \ln \f{\de}{2} } { \ga }, \; \; \mu
< \f{ \ga}{ \bs{n} } \ri \}$ is an empty set, thus the lemma is of course true.

Now we consider Case (ii).  By Lemma \ref{lem8}, there exists $\vep^* \in \li (0, \f{1}{\mu} - 1 \ri )$ such that {\small $\mscr{H} ( (1 +
\vep^*) \mu, \mu ) = \f{ \ln \f{\de}{2} } { \ga }$}.  By Lemma \ref{lem3},
\[
\Pr \li \{ \bs{n} \leq  \f{ \ga } { \mu (1 + \vep^*) }  \ri \} \leq \exp \li ( \ga  \mscr{H} ( (1 + \vep^*) \mu, \mu ) \ri ) = \f{\de}{2}.
\]
To show Lemma \ref{lem10}, it suffices to show
\[
\li \{  \mscr{H} \li ( \f{ \ga}{ \bs{n} },  \mu \ri ) \leq \f{ \ln \f{\de}{2} } {  \ga }, \; \; \mu < \f{ \ga}{ \bs{n} } \ri \} \subseteq \li \{
\bs{n} \leq  \f{ \ga } { \mu (1 + \vep^*) }  \ri \}.
\]
Let $\om \in \li \{  \mscr{H} \li ( \f{ \ga}{ \bs{n} },  \mu \ri ) \leq \f{ \ln \f{\de}{2} } {  \ga }, \; \; \mu < \f{ \ga}{ \bs{n} } \ri \}$
and $n = \bs{n} (\om)$. Let
\[
\vep =  \f{ \ga } { n \mu } - 1.
\]
Then,
\[
0 < \vep < \f{1}{\mu} - 1
\]
as a result of $n \geq \ga > n \mu$ and $0 < \mu < 1$.   By the definition of $\vep$, we have $\f{ \ga } { n} = (1 + \vep) \mu$ and thus
\[
\mscr{H} ( (1 + \vep) \mu, \mu ) \leq \f{ \ln \f{\de}{2} } {  \ga }.
\]
Hence, $\mscr{H} ( (1 + \vep) \mu, \mu ) \leq \mscr{H} ( (1 + \vep^*) \mu, \mu )$ and by Lemma \ref{lem2},
\[
\vep \geq \vep^*.
\]
By this inequality and  the definition of $\vep$,
\[
n = \f{ \ga } { \mu (1 + \vep) } \leq \f{ \ga } { \mu (1 + \vep^*) }.
\]
So,
\[
\om \in \li \{ \bs{n} \leq  \f{ \ga } { \mu (1 + \vep^*) }  \ri \}.
\]
This shows the inclusion relationship of the sets. The lemma is thus proved.
 \epf

 \bsk

 Now we are in a position to prove Theorem \ref{thm1}.  First, we shall  show that the confidence limits in Theorem \ref{thm1} are well-defined.  To this end, we
 need to show the existence and uniqueness of  $\udl{\bs{\mu}} (\om)$ and  $\ovl{\bs{\mu}} (\om)$  for any $\om \in \Om$.  For the lower
 confidence limit, the existence and uniqueness of  $\udl{\bs{\mu}} (\om)$ follows from Lemma \ref{lem6}, since $0 < \f{\ga}{ \bs{n} (\om) } \leq 1$ for any $\om \in \Om$.
For the upper confidence limit,  it is obvious that $\ovl{\bs{\mu}} (\om) = 1$ for $\om \in \Om$ such that $\bs{n} (\om) \leq \ga +
 1$.  For $\om \in \Om$ such that $\bs{n} (\om) > \ga +
 1$,  we have $0 < \f{ \ga } { \bs{n} (\om) - 1 } < 1$, and  the existence and uniqueness of $\ovl{\bs{\mu}} (\om)$ follows from Lemma
 \ref{lem6}.

Second,  we shall show that $\Pr \{ \mu \geq  \ovl{\bs{\mu}} \} \leq \f{\de}{2}$.
 By the definition of $\ovl{\bs{\mu}}$, we have
\bee \li \{ \mu \geq \ovl{\bs{\mu}} \ri \} & \subseteq & \li \{ \mscr{H} \li ( \f{ \ga}{ \bs{n} - 1 }, \ovl{\bs{\mu}} \ri ) = \f{ \ln \f{\de}{2}
} { \ga }, \qu \ovl{\bs{\mu}}
> \f{\ga}{ \bs{n} - 1 }, \qu \bs{n} > \ga + 1, \qu \mu \geq \ovl{\bs{\mu}} \ri \}\\
& \subseteq &  \li \{  \mscr{H} \li ( \f{ \ga}{ \bs{n} - 1},  \mu \ri ) \leq \f{ \ln \f{\de}{2} } {  \ga }, \; \; \mu \geq \f{ \ga}{ \bs{n} -
1}, \; \bs{n} > \ga + 1 \ri \} \eee where the second inclusion relationship can be shown as follows.

Let $\om \in \li \{ \mscr{H} \li ( \f{ \ga}{ \bs{n} - 1 }, \ovl{\bs{\mu}} \ri ) = \f{ \ln \f{\de}{2} } { \ga }, \qu \ovl{\bs{\mu}}
> \f{\ga}{ \bs{n} - 1 }, \qu \bs{n} > \ga + 1, \qu \mu \geq \ovl{\bs{\mu}} \ri \}$ and
\[
n = \bs{n} (\om), \qqu \ovl{\mu} = \ovl{\bs{\mu}} (\om).
\]
Then,  $\mu  \geq \ovl{\mu} > \f{\ga}{ n - 1 } > 0$ and $\mscr{H} \li ( \f{ \ga}{ n -1 },  \mu \ri ) \leq \mscr{H} \li ( \f{ \ga}{ n - 1 },
\ovl{\mu} \ri ) = \f{ \ln \f{\de}{2} } {  \ga }$ as a result of Lemma \ref{lem6}.  It follows that
\[
\om \in  \li \{  \mscr{H} \li ( \f{ \ga}{ \bs{n} - 1},  \mu \ri ) \leq \f{ \ln \f{\de}{2} } {  \ga }, \; \; \mu \geq \f{ \ga}{ \bs{n} - 1}, \;
\bs{n} > \ga + 1 \ri \}
\]
and the second inclusion relationship is true.   Applying Lemma \ref{lem9}, we have
\[
\Pr \li \{ \mu \geq \ovl{\bs{\mu}} \ri \} \leq \Pr \li \{  \mscr{H} \li ( \f{ \ga}{ \bs{n} - 1},  \mu \ri ) \leq \f{ \ln \f{\de}{2} } {  \ga },
\; \; \mu \geq \f{ \ga}{ \bs{n} - 1}, \; \bs{n} > \ga + 1 \ri \} \leq \f{\de}{2}.
\]

Third,  we shall show that  $\Pr \{ \mu \leq  \udl{\bs{\mu}} \} \leq \f{\de}{2} $.  By the definition of $\udl{\bs{\mu}}$, we have \bee \li \{
\mu \leq \udl{\bs{\mu}} \ri \} & \subseteq  &  \li \{ \mscr{H} \li ( \f{ \ga}{ \bs{n} }, \udl{\bs{\mu}} \ri ) = \f{ \ln \f{\de}{2} } {  \ga },
\qu
\udl{\bs{\mu}} < \f{\ga}{ \bs{n} }, \qu \mu \leq \udl{\bs{\mu}} \ri \}\\
& \subseteq & \li \{  \mscr{H} \li ( \f{ \ga}{ \bs{n} },  \mu \ri ) \leq \f{ \ln \f{\de}{2} } {  \ga }, \; \; \mu < \f{ \ga}{ \bs{n} } \ri \}
\eee where the second inclusion relationship can be shown as follows.

Let $\om \in \li \{ \mscr{H} \li ( \f{ \ga}{ \bs{n} }, \udl{\bs{\mu}} \ri ) = \f{ \ln \f{\de}{2} } {  \ga }, \qu \udl{\bs{\mu}} < \f{\ga}{
\bs{n} }, \qu \mu \leq \udl{\bs{\mu}} \ri \}$ and
\[
n = \bs{n} (\om), \qqu \udl{\mu} = \udl{\bs{\mu}} (\om).
\]
Then,  $\mu \leq \udl{\mu} < \f{\ga}{ n }$  and $\mscr{H} \li ( \f{ \ga}{n}, \mu \ri ) \leq \mscr{H} \li ( \f{ \ga}{ n }, \udl{\mu} \ri ) = \f{
\ln \f{\de}{2} } {  \ga }$ as a result of Lemma \ref{lem6}. Hence,
\[
\om  \in \li \{  \mscr{H} \li ( \f{ \ga}{ \bs{n} },  \mu \ri ) \leq \f{ \ln \f{\de}{2} } {  \ga }, \; \; \mu < \f{ \ga}{ \bs{n} } \ri \}
\]
and the second inclusion relationship is true.   Applying Lemma \ref{lem10}, we have
\[
\Pr \{ \mu \leq  \udl{\bs{\mu}} \} \leq  \Pr \li \{  \mscr{H} \li ( \f{ \ga}{ \bs{n} },  \mu \ri ) \leq \f{ \ln \f{\de}{2} } {  \ga }, \; \; \mu
< \f{ \ga}{ \bs{n} } \ri \} \leq \f{\de}{2}.
\]

Finally, Theorem \ref{thm1}  is justified by invoking the Bonferroni's inequality.

\sect{Proof Theorem 2}

To show Theorem \ref{thm2}, we need a modified version  of Lemma \ref{lem4} as follows.

\beL \la{lem11} For any $\ga > 0 $ and $\vep \in (0,1)$,
\[
\Pr \li \{  \bs{n} \geq \f{\ga}{\mu (1 - \vep)}  \ri \} \leq \exp \li ( \ga \mscr{H} \li ( \mu (1 -\vep), \mu \ri ) \ri ).
\]
\eeL

\bpf

Since $\bs{n}$ is an integer, we have
\[
\Pr \li \{  \bs{n} \geq \f{ \ga } { \mu (1 - \vep) } \ri \} =
 \Pr \li \{  \bs{n} \geq \li \lc \f{ \ga } { \mu (1
-\vep) } \ri \rc \ri \} =  \Pr \li \{  \bs{n} \geq  \f{ \ga } { \mu (1 - \ze^* ) }  \ri \}
\]
with $\ze^*$ satisfying {\small $\f{ \ga } { \mu (1 - \ze^* ) }  = \li \lc \f{ \ga } { \mu (1 - \vep ) } \ri \rc$}. Clearly, $1 > \ze^* \geq
\vep
> 0$.  Let {\small $m = \f{ \ga } { \mu (1 - \ze^* ) }$}. Then, $m$ is a positive integer and
\[
\Pr \{  \bs{n} \geq m \} = \Pr \{ X_1 + \cd + X_m \leq \ga \} = \Pr \{ \ovl{X} \leq z \}
\]
where {\small $\ovl{X} = \f{ \sum_{i=1}^m X_i } { m }$} and $z = (1 - \ze^*) \mu$.  Applying Lemma \ref{lem1}, we have {\small \bee \Pr \li \{
\bs{n} \geq \f{ \ga }
 { \mu (1 - \vep) }   \ri \}  & = & \Pr \{ \ovl{X} \leq z \} \leq  \exp ( m z
\mscr{H} ( z, \mu) )\\
&  = & \exp ( \ga \mscr{H} ( (1 - \ze^*) \mu, \mu) ). \eee} Note that $\mscr{H} ( (1 - \ze^*) \mu, \mu) \leq \mscr{H} ( (1 - \vep) \mu, \mu)$ as
a result of $1 > \ze^* \geq \vep > 0$ and Lemma \ref{lem2}. Hence,
\[
\Pr \li \{ \bs{n} \geq  \f{ \ga }
 { \mu (1 - \vep) }   \ri \} \leq \exp ( \ga \mscr{H} ( (1 - \vep) \mu, \mu) ).
\]

\epf

The remainder of the proof is similar to that of Theorem 1 and is thus omitted.

\end{document}